\documentclass[a4paper,12pt]{article}

\usepackage{amsmath,latexsym,amssymb}
\usepackage{bezier}
\usepackage[german,english]{babel}
\newtheorem{thm}{Theorem}[section]
\newtheorem{prop}[thm]{Proposition}
\newtheorem{lemma}[thm]{Lemma}

\newcommand{\proof}[1][]{{\it Proof#1: }}
\newcommand{\qed}[1][3mm]{\hspace*{\fill} $\Box$ \vspace{#1}}

\newcommand{\ON}{\operatorname}

\renewcommand{\b}{\beta}

\newcommand{\s}{\sigma}

\newcommand{\inj}{\hookrightarrow}
\newcommand{\surj}{\,\to\!\!\!\!\!\!\!\!\to}

\newcommand{\inv}{^{^{-1}}}
\newcommand{\invv}{^{-2}}

\newcommand{\rfami}{{\cal R}}

\newcommand{\cutoff}[1]{}
\newcommand{\labell}[1]{\label{#1}}

\newcommand{\br}{\ON{Br}}
\newcommand{\pbr}{\ON{PBr}}

\renewcommand{\iff}{\Leftrightarrow}

\newcommand{\ul}[1]{\underline{#1}}

\begin{document}

\title{Presentations of subgroups of the braid group
generated by powers of band generators}

\author{Michael L\"onne}

\maketitle

\begin{abstract}
According to the Tits conjecture proved by Crisp and Paris,
\cite{cp}, the subgroups of the braid group generated by
proper powers of the Artin elements $\s_i$ are 
presented by the commutators of generators which are
powers of commuting elements.
Hence they are naturally presented as right-angled Artin groups.

The case of subgroups generated by powers of the band generators
$a_{ij}$ is more involved. We show that the groups are
right-angled Artin groups again, if all generators are proper
powers with exponent at least $3$.
We also give a presentation in cases at the other extreme, when
all generators occur with exponent $1$ or $2$, which is
far from being that of a right-angled Artin group.
\end{abstract}

%
%
%
%
%
%

\section{Introduction}

Inspired by the Tits conjecture and its solution by Crisp and Paris, \cite{cp},
we investigate subgroups of the braid groups $\br_n$ generated by powers of the \emph{band generators}
$a_{ij}$.
They are the generators in the 'dual' or BKL presentation of $\br_n$:
$$
\left\langle 
a_{ij},\,1\leq i< j \leq n\,\left |
\begin{array}[c]{cl}
a_{ij}a_{kl}=a_{kl}a_{ij} & \mbox{\rm if } \ (k-i)(k-j)(l-i)(l-j)>0\\
a_{ij}a_{ik}=a_{jk}a_{ij}=a_{ik}a_{jk} & {\rm if\ }\ 1\leq i<j<k\leq n.
\end{array}
\right.\right\rangle 
$$
and can be identified with the following braid diagrams:

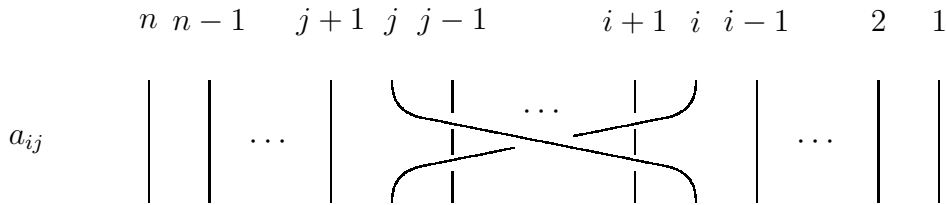
\unitlength=8mm
\begin{picture}(15,4.5)(-3,0)

\put(-1,2){\makebox(0,0){$a_{ij}$}}

\put(1,4){\makebox(0,0){$n$}}
\put(2,4){\makebox(0,0){$n-1$}}
\put(4,4){\makebox(0,0){$j+1$}}
\put(5,4){\makebox(0,0){$j$}}
\put(6,4){\makebox(0,0){$j-1$}}
\put(9,4){\makebox(0,0){$i+1$}}
\put(10,4){\makebox(0,0){$i$}}
\put(11,4){\makebox(0,0){$i-1$}}
\put(13,4){\makebox(0,0){$2$}}
\put(14,4){\makebox(0,0){$1$}}

\bezier{100}(1,1)(1,2)(1,3)
\bezier{100}(2,1)(2,2)(2,3)
\bezier{100}(4,1)(4,2)(4,3)
\bezier{100}(11,1)(11,2)(11,3)
\bezier{100}(13,1)(13,2)(13,3)
\bezier{100}(14,1)(14,2)(14,3)

\bezier{100}(5,1)(5,1.5)(5.5,1.6)
\bezier{100}(9.5,2.4)(10,2.5)(10,3)
\bezier{100}(5.5,1.6)(6,1.7)(7,1.9)
\bezier{100}(8,2.1)(9,2.3)(9.5,2.4)

\bezier{60}(5,3)(5,2.5)(5.5,2.4)
\bezier{190}(5.5,2.4)(7.5,2)(9.5,1.6)
\bezier{60}(9.5,1.6)(10,1.5)(10,1)

\bezier{50}(6,1)(6,1.2)(6,1.5)
\bezier{50}(6,1.8)(6,2)(6,2.2)
\bezier{50}(6,2.5)(6,2.8)(6,3)

\bezier{50}(9,1)(9,1.2)(9,1.5)
\bezier{50}(9,1.8)(9,2)(9,2.2)
\bezier{50}(9,2.5)(9,2.8)(9,3)

\put(3,2){\makebox(0,0){$\cdots$}}
\put(12,2){\makebox(0,0){$\cdots$}}
\put(7.5,2.5){\makebox(0,0){$\cdots$}}

\end{picture}
\\
So for a given symmetric matrix $M=(m_{ij})$ of exponents, the \emph{Coxeter datum},
we define
$$
E_{(m_{ij})}\,=\,\left\langle a_{ij}^{m_{ij}},\,1\leq i< j \leq n\right\rangle \,\subset\br_n
$$
As usual a vanishing exponent $m_{ij}=0$ yields the identity element,
thus $E_M$ is in fact generated by the powers of band generators $a_{ij}$ with
$m_{ij}\neq0$.

In that generality the groups $E_M$ were introduced by Kluitmann, cf. \cite{klui},
in conjunction with the more prominent Coxeter and Artin group associated with $M$:
\begin{eqnarray*}
C_M & \cong & \langle s_1,...,s_n\, | \, s_i^2, (s_is_j)^{m_{ij}}\,\rangle \\
A_M & \cong & \langle t_1,...,t_n\, | \, 
\underbrace{t_i t_j t_i...}_{m_{ij}\text{ factors}} = \underbrace{t_j t_i t_j...}_{m_{ij} \text{ factors}}\,\rangle
\end{eqnarray*}
\begin{quote}
(Note that our entries $m_{ij}=0$ in the Coxeter datum serve the same purpose as
entries $m_{ij}=\infty$, which are more commonly used.)
\end{quote}

For $M$ \emph{of large type}, ie.\ $m_{ij}\neq1,2$, Kluitmann proved
(and conjectured for \emph{non-redundant} $M$, ie.\ $m_{ij}\neq1$)
that $E_M$ is the stabiliser of $(s_1,\dots, s_n)$
for the natural action of $\br_n$
on $C_M\times\cdots\times C_M$.
If $m_{ij}=1$ is not excluded, $E_M$ may only be expected to be
the stabiliser of $(t_1,\dots,t_n)$ for the action on $A_M\times
\cdots\times A_M$, cf.\ \cite{hurwitz}.

For $M$ \emph{of ADE-type} D\"orner \cite{Doe} proved Kluitmanns conjecture.
The stabiliser group and hence $E_M$ has been shown by Looijenga, cf \cite{loo}
and others to be the fundamental group of the complement of the bifurcation set
in the truncated versal unfolding of the hypersurface singularity
of the same type (ADE).

In any case $E_M$ as a subgroup of a braid group $\br_n$ acts freely
and discontinously on a contractible CW-complex of dimension $n-1$,
and therefore deserves our special interest.
\\

While we would like to find a presentation of $E_M$ in terms of its
natural generators in general, we can provide two partial results in this
paper.

Our first result is to give the presentation of $E_M$ in case
of matrices $M$ of large type:

\begin{thm}
\labell{large}
In case of matrices $M$ of large type, 
ie. $m_{ij}\neq 1,2$ for all $1\leq i<j\leq n$, the
elements $b_{ij}=a_{ij}^{m_{ij}}$ with $m_{ij}\geq 3$ generate the subgroup $E_M$ 
with presentation
$$
\bigg\langle\, b_{ij},  m_{ij}\neq 0 \quad\left|\quad
b_{ij} b_{kl} = b_{kl} b_{ij},
\text{ if } a_{ij} a_{kl} = a_{kl} a_{ij}\,\bigg\rangle_.\right.
$$
\end{thm}

If $m_{ij}=0$ for all $i,j$ such that $|i-j|\geq2$ this covers the braid group case
of the theorem of Crisp and Paris.
\\

Second we determine a presentation in the cases at
the other extreme, when only exponents $1$ and $2$
are allowed. In this case several
distinct matrices may lead to the same subgroup $E_M$. Therefore
we include a condition on $M$ which give a distinguished matrix
in each case.

\begin{thm}
\label{manf}
Given a symmetric matrix $M$ with entries $m_{ij}\in\{1,2\}$
such that $m_{ik}$ is equal $1$ if $m_{ij}$ and $m_{jk}$ are,
a presentation of $E_M$ generated by $b_{ij}=a_{ij}^{m_{ij}}$
is obtained
imposing the following relations:\\[2mm]
\begin{tabular}{rcp{7.2cm}}
i) &
$
b_{il} b_{jk} = b_{jk} b_{il}
$ &
for $i<j<k<l$ or $j<k<i<l$\\[2mm]
ii) &
$
b_{jl} a_{kl}^2 b_{ik} a_{kl}\invv 
= a_{kl}^2 b_{ik} a_{kl}\invv b_{jl}
$  &
for $i<j<k<l$ and $a_{kl}^2=b_{kl}$ resp.\ $b_{kl}^2$, 
if $m_{kl}=2$ resp.\ $1$.\\[2mm]
iii) &
$
b_{ij} b_{ik} = b_{jk} b_{ij} ,b_{ij} b_{ik} b_{jk} = b_{ik} b_{jk} b_{ij}
$ &
for $i,j,k$ in cyclic order and $m_{ij}=1$, $m_{ik}=m_{jk}=2$,\\[2mm]
iv) &
$
b_{ij} b_{ik} b_{jk} = b_{jk} b_{ij} b_{ik} = b_{ik} b_{jk} b_{ij}
$ &
for $i<j<k$ and $m_{ij}=m_{ik}=m_{jk}=2$,\\[2mm]
v) &
$
b_{ij} b_{ik} = b_{jk} b_{ij} = b_{ik} b_{jk},
$ &
for $i<j<k$ and $m_{ij}=m_{ik}=m_{jk}=1$,
\end{tabular}
\end{thm}

Such groups $E_M$ can be identified with subgroups of $\br_n$
preserving a partition $P$ of $\{1,...,n\}$ via the natural surjection
$\br_n\to S_n$. These groups were considered before by Manfredini \cite{manf}
who gave presentations which are more economical in terms of
generators than ours. In contrast it is much easier to state and prove the
relations of our presentations.

In particular the pure braid group $\pbr_n$ corresponds to
the matrix $M$ with all entries $m_{ij}=2$. We thus get
a presentation of $\pbr_n$, which differs from the 
classical presentation of Artin \cite{ar}
but recovers in fact the earlier one by Burau \cite{bu}.
\\

Of the following three sections we devote the first two to the
proofs of the two theorems, while in the last sections we want
to suggest an approach in the case with $M$ of ADE-type.

\section{The case of matrices of large type}

In their approach Crisp and Paris exploit
a geometric action of Artin groups on a
fundamental groupoid and the solution of 
the word problem for right-angled Artin groups.

We follow their strategy but in place of their action
we pick a natural geometric representation of $\br_n$, 
the identification with the mapping class group
of the punctured disc.
Then we get an induced action on the fundamental group of the punctured disc,
which we identify with the free group $F_n=\langle t_1,...,t_n\,|\,\quad\rangle$.
The right Hurwitz actions of $\br_n$ on $F_n$ and the universal Coxeter quotient
$\tilde C_n=\langle s_1,...,s_n\,|\,s_i^2=1\,\rangle$
are given by
$$
\begin{matrix}
\begin{array}{rcl}
(t_i) \s_j & = &
\left\{
\begin{array}{ll}
t_{i-1} & \text{if } i=j+1\\
t_it_{i+1} t_i\inv & \text{if } i=j\\
t_i & \text{else}
\end{array}
\right.
\end{array}
&
\begin{array}{rcl}
(s_i) \s_j & = &
\left\{
\begin{array}{ll}
s_{i-1} & \text{if } i=j+1\\
s_i s_{i+1} s_i & \text{if } i=j\\
s_i & \text{else}
\end{array}
\right.
\end{array}
\end{matrix}
$$

We must now define the terms to describe some decisive features of the action of generators
$a_{jk}^m$ on elements in $\tilde C_n$, which we represent by words without repetition in the alphabet $s_1,...,s_n$.
Occasionally a transformation rule will result in more general words, but reduction, the process of removing
a pair of adjacent identical letters, will unambiguously stop at the
representing word without repetitions.

\paragraph*{Definition}
A subword $w_0$ of a word $w$ in the alphabet $s_1,...,s_n$ 
is called a \emph{$jk$-subword}
if it consists of letters $s_j$ and $s_k$ only, it is called \emph{long} if
it contains at least four letters.

\paragraph{Definition}
The \emph{$jk$-factorisation} of a word $w$ is given by the unique factorisation
$$
w\quad=\quad
w_0s_{i_1}w_1\cdots s_{i_\ell} w_\ell,
$$
such that each $w_\nu$ is a $jk$-subword of $w$, possibly of length $0$,
and such that $1\leq i_\nu \leq n$, $i_\nu\neq j,k$ for all $\nu$.

\paragraph*{Definition}
A $jk$-subword $w_\nu$ is called \emph{critical} if either of the following condition holds
\begin{enumerate}
\item
the length $|w_\nu|$ is odd and $w_\nu$ has two identical adjacent letters, distinct from $s_j,s_k$,
\item
the length $|w_\nu|$ is even, possibly $0$, and $w_\nu$ has adjacent letters $s_i$ and $s_l$
with $il$ and $jk$ crossing, that is $j,k$ are either in the even
or the odd positions for the natural order of the four numbers $i,j,k,l$.
\end{enumerate}

Note that being critical is \emph{not} an absolute notion but depends on the
'embedding' as a subword, more precisely on the adjacent letters.

\begin{prop}
\labell{trans}
Suppose that a word $w$ has a $jk$-factorisation  
$w_0s_{i_1}w_1\cdots s_{i_\ell} w_\ell$ with no
$jk$-subword $w_\nu$ of length $2|m|$.
Then the $jk$-factorisation of $(w)a_{jk}^m$ can be written as
$$
w_0's_{i_1}w_1'\cdots s_{i_\ell} w_\ell',
$$
where $w_\nu'$ is critical with $|w_\nu'|+|w_\nu|\geq 2|m|$ if $w_\nu$ is critical.
\end{prop}

\proof
From the Hurwitz action we deduce the following action of a band generator and its powers on any letter
using the identity $(s_js_k)=(s_ks_j)\inv$:
\begin{eqnarray*}
(s_i) a_{jk}^m  & = &
\left\{
\begin{array}{ll}
s_{i}, &  i<j \text{ or } i>k\\
s_j(s_j s_k)^{-m}= (s_j s_k)^m s_j, &  i=j\\
(s_j s_k)^m s_i (s_j s_k)^{-m}, & j<i<k\\
(s_j s_k)^m s_k = s_k (s_j s_k)^{-m}, & i=k 
\end{array}
\right.\\
\end{eqnarray*}
They imply the following transformation of various $jk$-subwords, $p\geq0$:
\begin{eqnarray*}
((s_j s_k)^p)a_{jk}^m & = & (s_j s_k)^p\\
((s_k s_j)^p)a_{jk}^m & = & (s_j s_k)^{-p}\\
((s_j s_k)^p  s_j)a_{jk}^m & = & (s_j s_k)^p(s_j s_k)^m s_j \quad=\quad (s_j s_k)^{p+m} s_j\\
((s_k s_j)^p  s_k)a_{jk}^m & = & (s_j s_k)^{-p}(s_j s_k)^m s_k \quad=\quad (s_j s_k)^{-p+m} s_k
\end{eqnarray*}
Hence we obtain the word $w'=(w)a_{jk}^m$ by first applying these transformation to the factors $w_\nu,s_{i_\nu}$
and then removing adjacent pairs of identical letters.

By the first step we get in general a \emph{non-reduced} word.
Its $jk$-factorisation has the same number of terms as that of $w$,
since the number of letters distinct from $s_j,s_k$ is preserved.

In the next step we want to remove all repetitious pairs in the $jk$-subwords.
Each $jk$-subword between $s_{i_\nu}$ and $s_{i_{\nu+1}}$
consist of contributions from three sources:
\begin{enumerate}
\item
the first from $s_{i_\nu}$: $(s_js_k)^{-m}$ if
$j<i_\nu<k$, nothing otherwise.
\item
the second from $w_\nu$: $w_\nu$ if $|w_\nu|$ is even, $(s_j s_k)^m w_\nu$ if $|w_\nu|$ is odd.
\item
the third from $s_{i_{\nu+1}}$: $(s_js_k)^{m}$ if $j<i_{\nu+1}<k$, nothing otherwise.
\end{enumerate}

So the $jk$-subwords $w_\nu'$ may be written as either $(s_js_k)^m w_\nu$, $(s_j s_k)^{-m}w_\nu$
or $w_\nu$.
Reducing further we get eventually the $jk$-word $w_\nu'$, which may be empty only if $w_\nu$ is
since $|w_\nu|\neq 2|m|$ is assumed.

But then also the word $w_0's_{i_1}w_1'\cdots s_{i_\ell} w_\ell'$ 
has no repetitions, 
thus our claim on the $jk$-factorisation of $w'$ is proved.

We are left to consider the two particular cases when $w_\nu$ is critical.
If $i_\nu=i_{\nu+1}$ and $|w_\nu|$ is odd, then $w_\nu'=(s_j s_k)^{\pm m}w_\nu$,
hence $w_\nu'$ is critical and $|w_\nu| + |w_\nu'| \geq 2|m|$.
The same conclusion holds if
$i_\nu i_{\nu+1}$ and $jk$ are crossing and $|w_\nu|$ is even.
\qed

\begin{prop}
\labell{7}
Suppose that $w$ is without long $il$-subword
while $(w)a_{il}^m$, with $|m|\geq 3$ has a long $jk$-subword.
Then one of the following two conditions holds:
\begin{enumerate}
\item
$il$ and $jk$ are identical,
\item
$il$ and $jk$ are non-crossing and $w$ has a long $jk$-subword.
\end{enumerate}
\end{prop}

\proof
By the previous proposition we may write the $il$-factorisations of
$w$ and its image $w'=(w)a_{il}^m$ as
$$
w_0s_{i_1}w_1\cdots s_{i_\ell} w_\ell
\quad\text{resp.}\quad
w_0's_{i_1}w_1'\cdots s_{i_\ell} w_\ell'.
$$
Any long $jk$-subword in $w'$ is either
\begin{enumerate}
\item
a long $il$-word, if $il$ and $jk$ are identical,
\item
an alternating sequence of at least two $il$-words of length $1$
and two copies of a letter not in $\{s_i,s_l\}$, 
if $il$ and $jk$ are neither identical nor disjoint,
\item
an alternating sequence of letters $s_j,s_k$ interspersed by 
$il$-subwords of length $0$,
if $il$ and $jk$ are disjoint.
\end{enumerate}
The first case is accounted for by the claim.

In the second case $w'$ contains a maximal $il$-subword $w_\nu'$ 
of length $1$ with identical adjacent letters.
Hence $w_\nu'$ is critical. But then by prop. \ref{trans} also $w_\nu$ must be critical and of length $|w_\nu|\geq2|m|-|w_\nu'|\geq 5$
in violation of our assumptions.

In the last case each interspersed $il$-subword $w_\nu'$ may not be critical, otherwise it would originate in a critical $il$-subword
$w_\nu$ of length $|w_\nu|\geq2|m|-|w_\nu'|\geq 6$
which again is excluded by our assumption.
Since the adjacent letters are $s_j$ and $s_k$, we conclude
that $il$ and $jk$ are non-crossing, as claimed.
\qed

For the ensuing argument we introduce some notation to handle
elements of the right-angled Artin groups $G_M$ defined on generators
in bijection to elements of $T_M$:
\begin{eqnarray*}
G_M & = &
\left\langle\, b_\tau,\, \tau \in T_M\:|\:
b_\tau b_\s = b_\s b_\tau \text{ if $\tau,\s$ are non-crossing}\;\right\rangle\\
T_M & = &
\left\lbrace \,ij\:|\: 1\leq i<j\leq n,\,m_{ij}\neq0\, \right\rbrace 
\end{eqnarray*}

An \emph{expression} in the letters $b_{\s}$ is a sequence with non-vanishing exponents $p_i$
$$
W=(b_{\s_1}^{p_1},b_{\s_2}^{p_2},..., b_{\s_\ell}^{p_\ell})
$$
The index $\ell=\ell(W)$ is called the length of the
expression $W$.
Given $\b$ in $E_M$ we call $W$ an expression for $\b$ if
$$
\b\quad=\quad
a_{\s_1}^{p_{1}m_{\s_{1}}}\cdots
a_{\s_\ell}^{p_{\ell}m_{\s_{\ell}}}
$$
As in \cite{cp} the following terminology is based on Brown.
Consider an expression $W$ as above. Suppose that there exists $i\in \{1,...,\ell\}$ such that
$b_{\s_i}=b_{\s_{i+1}}$, put 
$$
W'\quad=\quad
\left\lbrace 
\begin{array}{ll}
(b_{\s_1}^{p_1},b_{\s_2}^{p_2},...,b_{\s_{i-1}}^{p_{i-1}},b_{\s_i}^{p_i+p_{i+1}},b_{\s_{i+2}}^{p_{i+2}},...,b_{\s_\ell}^{p_\ell}) &
\text{if } p_i+p_{i+1}\neq 0,\\[2mm]
(b_{\s_1}^{p_1},b_{\s_2}^{p_2},...,b_{\s_{i-1}}^{p_{i-1}},b_{\s_{i+2}}^{p_{i+2}},...,b_{\s_\ell}^{p_\ell}) &
\text{if } p_i+p_{i+1}=0.
\end{array}
\right.
$$
We say that $W'$ is obtained from $W$ via an elementary operation of type $I$.
This operation shortens the length of an expression by $1$ or $2$.

Suppose that there exists $i\in\{1,...,\ell-1\}$ such that $\s_i,\s_{i+1}$ are non-crossing.
Put
$$
W''\quad=\quad
(b_{\s_1}^{p_1},b_{\s_2}^{p_2},...,b_{\s_{i-1}}^{p_{i-1}},b_{\s_{i+1}}^{p_{i+1}},b_{\s_i}^{p_i},
b_{\s_{i+2}}^{p_{i+2}},...,b_{\s_\ell}^{p_\ell}).
$$
We say that $W''$ is obtained from $W$ by an elementary operation of type $II$.
This operation leaves the length of an expression unchanged.

We shall say that $W$ is $M$-reduced if the length of $W$ can not be reduced by applying a sequence of elementary
operations. Clearly every element of the right-angled presented group has a reduced expression.

A reduced expression $W$ is said to end in $\tau \in T_M$ if it is related by a
sequence of operations of type $II$ to an expression 
$(b_{\s_1}^{p_1},b_{\s_2}^{p_2},...,b_{\s_\ell}^{p_\ell})$ in which $\s_\ell=\tau$.

\begin{prop}
\labell{9}
Suppose $X=(x_1,...,x_\ell), x_\nu=b_{\s_\nu}^{p_\nu},$ is a $M$-reduced expression for $\b\in E_M\subset\br_n$.
If $(s_i)\b$ contains a long $jk$-subword,
then $X$ ends in $\tau=jk$.
\end{prop}

\proof
We argue by induction on the length $\ell$ of $X$;
obviously for trivial $X$ there is nothing to be proved.

Hence given an expression $X$ of length $\ell>0$ we may assume that the claim holds
for $X'=(x_1,...,x_{\ell-1})$, which is a $M$-reduced expression for $\b'$,
where $\b=\b' a_{\s_{\ell}}^{mp_{\ell}}, m=m_{\s_{\ell}}$.

Note that $(s_i)\b'$ has no long $\s_\ell$-subword,
for otherwise $X'$ ends with $\s_\ell$ by induction, contrary to our assumption that $X$ is $M$-reduced.

Since $(s_i)\b=\left( (s_i)\b'\right)a_{\s_{\ell}}^{mp_{\ell}} $ contains a long $\tau$-subword,
we conclude with
prop. \ref{7} that either $\tau=\s_\ell$
or $\tau,\s_\ell$ are non-crossing and $(s_i)\b'$ contains a long $\tau$-subword.

In the first case the claim is obviously true.
In the second case, by induction, $X'$ ends in $\tau$, hence there is an
expression $(x_1',...,x_{\ell-1}')$ of $\b'$ with $\s_{\ell-1}'=\tau$, which is
obtained by operations of type $II$ from $X'$.
Hence $X$ transforms into $(x_1',...,x_{\ell-1}',x_{\ell})$.

Since $\tau,\tau_\ell$ are non-crossing
 a further operation of type $II$ yields
$$
(x_1',..., x_{\ell-2}',x_\ell^{}, x_{\ell-1}')
$$
which shows that $X$ ends in $\tau$ as claimed.
\qed

\begin{prop}
\labell{10}
Suppose $X=(x_1,...,x_\ell)$ is a non-trivial $M$-reduced expression for $\b$
with $\s_\ell={jk}$, $x_l=b_{jk}^{p_\ell}$.
Then $(s_j)\b\neq s_j$, in particular $G_M$ injects into $\br_n$.
\end{prop}

\proof
Again $X'=(x_1,...,x_{\ell-1})$ is an expression for $\b'=\b a_{jk}^{-p_\ell m_{jk}}$.
Suppose contrary to our claim that $(s_j)\b=s_j$, so
$$
(s_j)\b'=(s_j)\b a_{jk}^{-p_\ell m_{jk}}=(s_j) a_{jk}^{-p_\ell m_{jk}}.
$$
By the transformation rules in proposition \ref{trans}
the word on the right contains a long $jk$-subword, hence
the same is true on the left hand side.
With proposition \ref{9} we conclude that $X'$ ends in $\s_\ell={jk}$.
Now we have reached a contradiction to the hypothesis that $X$ is $M$-reduced,
thus our claim holds true.
\qed

The observation of the proposition concludes the proof of theorem \ref{large}
since $G_M$ now maps isomorphically onto $E_M$ identifying the respective
sets of generators.
\\

Our argument give also a new prove of the Tits conjecture in case
of the braid group and large exponents, ie.\ at least $3$.

\section{partition preserving braid subgroups}

We want to study the groups $E_M$ of theorem \ref{manf}
in terms of corresponding partitions $P$ of the set
$\{1,...,n\}$. Given such a partition we define its stabiliser
in the symmetric group $S_n$ to be
$$
S_{n,P}\quad =\quad \langle (ij)| i\sim_P j\rangle\:\subset\: S_n,
$$
where $\sim_P$ is the unique equivalence relation with set of
equivalence classes equal to $P$.

The significance of these notions for our problem lies in
the following correspondence:

\begin{prop}
\label{partit}
There is a natural bijection between matrices $M$
as in theorem \ref{manf} and the
partitions $P$ of $\{1,...,n\}$ such that there is a natural
diagram with exact rows:
$$
\begin{array}{ccccccccc}
1 & \to & \pbr_n & \to & E_M & \to & S_{n,P} & \to & 1\\
&& \| && \cap && \cap\\
1 & \to & \pbr_n & \to & \br_n & \to & S_{n} & \to & 1
\end{array}
$$
\end{prop}

\proof
To a partition $P$ we associate the unique matrix $M$
with $m_{ij}=1 \iff i\sim_P j$ and $m_{ij}=2$ otherwise.
The additional condition, $m_{ij}=m_{jk}=1\Rightarrow
m_{ik}=1$, holds for $M$ since $\sim_P$
is transitive. Reversely with $M$ we may define the relation
$i\sim_M j:\iff m_{ij}=1\vee i=j$, which obviously is reflexive and
symmetric, but also transitive due to the additional condition on $M$.
The set of equivalence classes yields the partition $P$ associated with $M$.

For the second claim we first convince ourselves that $S_{n,P}$ is
generated by all transpositions $(ij)$ with $i\sim_P j$.
For the associated $M$ obviously $E_M$ contains $\pbr_n$.
Moreover $(ij)$ has a preimage $b_{ij}$ in $E_M$ if and only if
$m_{ij}=1$, hence $E_M$ surjects onto $S_{n,P}$.
\qed

In order to prove theorem \ref{manf}
we want to set up an induction over $n$ which is well-based, since
the cases $n=1,n=2$ are obviously true.
Given a partition $P$ of $\{1,...,n\}$ we define the induced
partition $P'$ of $\{1,...,n-1\}$ and we denote by
$P'n$ the partition of $\{1,...,n\}$ for which $\{n\}$ is a part on
its own and all other part are those of $P'$.

Moreover we use $a_i$ as a shorthand for $a_{i\,n}$ 

\begin{lemma}
\label{combing}
Suppose $E_{P'}$ is presented as
$\langle b_{jk}, 1\leq j<k< n \,|\,\rfami'\rangle$,
then $E_{P'n}$ has a presentation by
generating elements 
$
a_{i}^2, b_{jk}
$
subject to relations in $\rfami'$ and \\[2mm]
\begin{tabular}{rcp{4.7cm}}
i) &
$
a^2_{i} b_{jk} = b_{jk} a^2_{i}
$ &
for $i<j<k< n$ or\hfill\phantom{x} $j<k<i< n$\\[2mm]
ii) &
$
a^2_{j} a^2_{k} b_{ik} a\invv_{k} = 
a^2_{k} b_{ik} a_{k}\invv a^2_{j}.
$ &
for $i<j<k< n$\\[2mm]
iii) &
$
b_{ij} a^2_{i} = a^2_{j} b_{ij},
b_{ij} a^2_{i} a^2_{j} = a^2_{i} a^2_{j} b_{ij}
$ &
for $i<j< n$ and $m_{ij}=1$,\\[2mm]
iv) &
$
a^2_{j} b_{ij} a^2_{i } = 
b_{ij} a^2_{i} a^2_{j} = a^2_{i} a^2_{j} b_{ij}
$ &
for $i<j< n$ and $m_{ij}=2$,
\end{tabular}\\[2mm]
In particular
the theorem \ref{manf} holds true for $E_{P'n}$,
if it holds true for $E_{P'}$.
\end{lemma}

\proof
All given relations are shown to hold by a straightforward calculation.
On the other hand we can argue with the following natural diagram
of split exact rows:
$$
\begin{array}{ccccccccc}
1 & \to & \langle a_{i\,n }^2\rangle &
\to & E_{ P'n} & \to & E_{P'} & \to & 1\\
& & \| & & \cap & & \cap \\
1 & \to & \langle a_{i\,n }^2\rangle &
\to & \br_{n-1,1} & \to & \br_{n-1} & \to & 1\\
\end{array}
$$
Hence it suffices to show that all relations obtained from the action
of the extension can be deduced from those given in the claim.
These relations are easily obtained by 'combing' the braid obtained
by conjugation.
The following list -- with $i<j<k< n$ and $m_{ij}=1$ in case $3,4$,
resp.\ $m_{ij}=2$ in case $5,6$, $m_{ik}=1$ in case $7$ and
$m_{ik}=2$ in case $8$ -- is exhaustive.

$$
\begin{array}{crclc}
(1)& b_{jk} a^2_{i\,n }b_{jk}\inv 
& = & 
a^2_{i\,n } \\[1mm]
&\iff & \multicolumn{2}{l}{ b_{jk} a^2_{i\,n } 
\quad = \quad 
a^2_{i\,n }b_{jk}} \\[1mm]
(2)& b_{ij} a^2_{k\,n } b_{ij}\inv
& = & 
a^2_{k\,n } \\[1mm]
&\iff & \multicolumn{2}{l}{b_{ij} a^2_{k\,n }  
\quad = \quad 
a^2_{k\,n } b_{ij}} \\[1mm]
(3)& b_{ij} a^2_{i\,n }b_{ij}\inv  
& = & 
a^2_{j\,n } 
& (m_{ij}=1)\\[1mm]
&\iff & \multicolumn{2}{l}{ b_{ij} a^2_{i\,n }
\quad = \quad 
a^2_{j\,n } b_{ij} } 
& \text{\it iii)}\\[1mm]
(4)& b_{ij}a^2_{j\,n } b_{ij}\inv  
& = &
a_{j\,n }\invv a^2_{i\,n } a^2_{j\,n } 
& (m_{ij}=1)\\[1mm]
&\iff & \multicolumn{2}{l}{ \ul{a^2_{j\,n }b_{ij}} a^2_{j\,n } 
\quad = \quad 
a^2_{i\,n }a^2_{j\,n }b_{ij} } 
& \text{\it iii)}\\[1mm]
&\iff & \multicolumn{2}{l}{  b_{ij} a^2_{i\,n }a^2_{j\,n }
\quad = \quad 
a^2_{i\,n }a^2_{j\,n } b_{ij}} 
& \text{\it iii)}\\[1mm]
(5)& b_{ij} a^2_{i\,n }b_{ij}\inv 
& = & 
a_{j\,n }\invv a^2_{i\,n } a^2_{j\,n }
& (m_{ij}=2)\\[1mm]
&\iff & \multicolumn{2}{l}{ a^2_{j\,n } b_{ij} a^2_{i\,n }
\quad = \quad 
a^2_{i\,n }a^2_{j\,n }b_{ij} } 
& \text{\it iv)}\\[1mm]
(6)& b_{ij} a^2_{j\,n } b_{ij}\inv 
& = & 
a_{j\,n }\invv a_{i\,n }\invv a^2_{j\,n }a^2_{i\,n } a^2_{j\,n } 
& (m_{ij}=2)\\[1mm]
&\iff & \multicolumn{2}{l}{ \ul{a^2_{i\,n }a^2_{j\,n } b_{ij}} a^2_{j\,n } 
\quad = \quad 
a^2_{j\,n }\ul{a^2_{i\,n } a^2_{j\,n }b_{ij}} } 
& \text{\it iv)}\\[1mm]
&\iff & \multicolumn{2}{l}{ a^2_{j\,n } b_{ij} a^2_{i\,n }a^2_{j\,n }
\quad = \quad 
a^2_{j\,n } b_{ij} a^2_{i\,n }a^2_{j\,n }} \\[1mm]
(7)& b_{ik} a^2_{j\,n } b_{ik}\inv 
& = &
 a_{k\,n }\invv a^2_{i\,n } a^2_{j\,n }a_{i\,n }\invv a^2_{k\,n }
& (m_{ij}=1)\\[1mm]
&\iff & \multicolumn{2}{l}{ a_{i\,n}\invv \ul{a^2_{k\,n } b_{ik}} a^2_{j\,n }  
\quad = \quad
 a^2_{j\,n }a_{i\,n }\invv \ul{a^2_{k\,n } b_{ik}} } 
& \text{\it iii)}\\[1mm]
&\iff & \multicolumn{2}{l}{ a^2_{j\,n } a^2_{k\,n } b_{ik} a_{k\,n }\invv 
\quad = \quad 
a^2_{k\,n } b_{ik} a_{k\,n }\invv a^2_{j\,n }} 
& \text{\it ii)}\\[1mm]
(8)& b_{ik}\inv a^2_{j\,n } b_{ik} 
& = & 
a^2_{i\,n }a^2_{k\,n } a_{i\,n }\invv a_{k\,n }\invv a^2_{j\,n } 
a^2_{k\,n } a^2_{i\,n }a_{k\,n }\invv a_{i\,n }\invv 
& (m_{ij}=2)\\[1mm]
&\iff & \multicolumn{2}{l}{ a^2_{j\,n } \ul{b_{ik} a^2_{i\,n }a^2_{k\,n }} a_{i\,n }\invv
a_{k\,n }\invv
\quad = \quad 
\ul{b_{ik} a^2_{i\,n }a^2_{k\,n }} a_{i\,n }\invv 
a_{k\,n }\invv a^2_{j\,n }} 
& \text{\it iv)}\\[1mm]
&\iff & \multicolumn{2}{l}{ a^2_{j\,n } a^2_{k\,n } b_{ik} (a^2_{i\,n }a_{i\,n }\invv) 
a_{k\,n }\invv
\quad = \quad 
a^2_{k\,n } b_{ik} (a^2_{i\,n }a_{i\,n }\invv) a_{k\,n }\invv 
a^2_{j\,n }}
& \text{\it ii)}
\end{array}
$$
We may now apply the claim of the main part of the lemma to the
particular case in which $\rfami'$ is the set of relations associated to
$E_{P'}$ and the generators $b_{jk}, 1\leq j < k <n$ in the claim of
the theorem and $\langle b_{jk}, j,k<n|\rfami'\rangle$
is a presentation of $E_{P'}$.

Then of course $E_{P'n}$ is presented by generators
$b_{ij}, i,j<n$ and $b_{i\,n}=a_{i\,n}^2$ subject to the relations
$\rfami'$ and those of the list, we just derived.
Since they correspond bijectively to relations of the theorem \ref{manf}
in case $P=P'n$ we have concluded our proof.
\qed

Now given a general $P$ we consider the group $G_M$ 
given by the presentation of theorem \ref{manf}
in terms of the matrix $M$ associated to $P$ by \ref{partit},
but with elements $\tilde b_{ij}$ to mark
the distinction with the subgroups of $\br_n$.
Let $I\subset\{1,...,n\}$ 
be the equivalence class of $n$ under $\sim_P$.
Then we denote by $H_M$ the subgroup of $G_M$
generated by elements $\{\tilde b_{ij},i<j<n\}$,
$\{\tilde b_i:=\tilde b_{in}, i\not\in I\}$
and $\{\tilde b_i^2:=\tilde b_{in}^2, i\in I\}$,
where we take $\tilde b_n=1$.

\begin{lemma}
\label{cosets}
$G_M$ can be given as a finite union of
right cosets of $H_M$ (with $\tilde b_n=1$):
$$
G_M \quad =\quad \bigcup_{t\in I}\, \tilde b_t H_M.
$$
\end{lemma}

\proof
Since each generator is either in
$H_M$ or in $\{\tilde b_t|t\in I\}$, all 
generators of $G_M$ belong to this union of
cosets and multiplication by any generator on the left
maps $H_M$ to one of the given cosets.
Hence it remains to prove that multiplication on the left
by any generator $\tilde b$ of $G_M$ maps an element
of any
of the non-trivial given cosets of $H_M$ into one
of the given cosets.
We refer by roman numbers to the relations of the theorem
which are used and freely move factors from and into $H_M$:
\begin{enumerate}
\item
$\tilde b= \tilde b_i, i\in I$
\begin{enumerate}
\item
$i<t$\\[-5.5mm]
$$
\tilde b_i \tilde b_t \stackrel{v)}{=} \tilde b_t \tilde b_{it} \in \tilde b_t H_M
$$
\item
$i=t$\\[-5.5mm]
$$
\tilde b_i \tilde b_t \stackrel{}{=} \tilde b_t^2 \in  H_M
$$
\item
$i>t$\\[-5.5mm]
$$
\tilde b_i \tilde b_t \in \tilde b_i \tilde b_t \tilde b_t\invv H_M
\stackrel{}{=} 
\tilde b_i \tilde b_t\inv H_M \stackrel{v)}{=}
\tilde b_t\inv \tilde b_{ti} H_M \stackrel{}{=}
\tilde b_t \tilde b_t\invv H_M
\in \tilde b_t  H_M
$$
\end{enumerate}
\item $\tilde b = \tilde b_i, i\not\in I$
\begin{enumerate}
\item
$i<t$\\[-5.5mm]
$$
\tilde b_i \tilde b_t \stackrel{iii)}{=} 
\tilde b_t \tilde b_{it} \in \tilde b_t  H_M
$$
\item
$i=t$ not possible since $t\in I, i\not\in I$,
\item
$i>t$\\[-5.5mm]
$$
\tilde b_i \tilde b_t \in
\tilde b_i \tilde b_t \tilde b_i H_M  \stackrel{iii)}{=} 
\tilde b_i \tilde b_{ti} \tilde b_t H_M \stackrel{iii)}{=}
\tilde b_t \tilde b_i \tilde b_{ti} H_M \stackrel{}{=}
 \tilde b_t H_M
$$
\end{enumerate}
\item
$\tilde b=\tilde b_{ij}$
\begin{enumerate}
\item $t<i<j$ or $i<j<t$
$$
\tilde b_{ij} \tilde b_t \stackrel{i)}{=} 
\tilde b_t \tilde b_{ij} \in \tilde b_t H_M
$$
\item $t=i<j, j\in I$, ie. $m_{ij}=1$
$$
\tilde b_{tj} \tilde b_t \stackrel{v)}{=} 
\tilde b_j \tilde b_{tj} \in \tilde b_j H_M
$$
\item $i<j=t, i\in I$, ie. $m_{ij}=1$
$$
\tilde b_{it} \tilde b_t \in
\tilde b_{it} \tilde b_t \tilde b_{it} H_M \stackrel{v)}{=} 
\tilde b_{it} \tilde b_i \tilde b_{t} H_M \stackrel{v)}{=} 
\tilde b_i \tilde b_{t}^2  H_M \stackrel{}{=} 
\tilde b_i  H_M
$$
\item $t=i<j, j\not\in I$, ie. $m_{ij}=2$
$$
\tilde b_{tj} \tilde b_t \stackrel{iii)}{=}
\tilde b_t \tilde b_j \in \tilde b_t H_M
$$
\item $i<j=t, i\not\in I$, ie. $m_{ij}=2$
$$
\tilde b_{it} \tilde b_t \in
\tilde b_{it} \tilde b_t \tilde b_{it} H_M \stackrel{iii)}{=} 
\tilde b_{it} \tilde b_i \tilde b_t H_M \stackrel{iii)}{=}
\tilde b_t \tilde b_{it} \tilde b_i H_M 
\stackrel{}{=} \tilde b_j  H_M
$$
\end{enumerate}
\item
$\tilde b= \tilde b_{ik}$, $i<t<k$,
\begin{enumerate}
\item $i,t\in I$, $k\not\in I$, i.e. $m_{ik}=2$, then with
$*)$: cases $(a)$ of $i)$ and $(c)$ of $ii)$ above
$$
\tilde b_{ik} \tilde b_t \stackrel{*)}{\in}
\tilde b_{ik} \tilde b_i^2 \tilde b_k \tilde b_t H_M \stackrel{iii)}{=}
\tilde b_i \tilde b_k \tilde b_{ik} \tilde b_i \tilde b_t H_M \stackrel{iii)}{=}
\tilde b_k \tilde b_{ik} \tilde b_i^2 \tilde b_t H_M \stackrel{*)}{=}
\tilde b_k \tilde b_{ik} \tilde b_k\inv \tilde b_t H_M
$$
$$
\stackrel{ii)}{=}
\tilde b_t \tilde b_k \tilde b_{ik} \tilde b_k\inv H_M \stackrel{}{=}
\tilde b_t H_M \hspace*{6cm}
$$
\item $t\in I$, $i,k\not\in I$, $m_{ik}=2$, then with
$*)$: cases $(a)$ and $(c)$ of $ii)$ above
$$
\tilde b_{ik} \tilde b_t \stackrel{*)}{\in}
\tilde b_{ik} \tilde b_i \tilde b_k \tilde b_t H_M \stackrel{iv)}{=}
\tilde b_k \tilde b_{ik} \tilde b_i \tilde b_t H_M \stackrel{*)}{=}
\tilde b_k \tilde b_{ik} \tilde b_k\inv \tilde b_t H_M \stackrel{ii)}{=}
\tilde b_t \tilde b_k \tilde b_{ik} \tilde b_k\inv H_M \stackrel{}{=}
\tilde b_t H_M
$$
\item $k,t\in I$, $i\not\in I$, i.e. $m_{ik}=2$, then with
$*)$: cases $(c)$ of $i)$ and $(a)$ of $ii)$ above
$$
\tilde b_{ik} \tilde b_t \stackrel{*)}{\in}
\tilde b_{ik} \tilde b_i \tilde b_k^2 \tilde b_t H_M \stackrel{iii)}{=}
\tilde b_k^2 \tilde b_{ik} \tilde b_i \tilde b_t H_M \stackrel{*)}{=}
\tilde b_k^2 \tilde b_{ik} \tilde b_k\invv \tilde b_t H_M \stackrel{ii)}{=}
\tilde b_t \tilde b_k^2 \tilde b_{ik} \tilde b_k\invv H_M \stackrel{}{=}
\tilde b_t H_M
$$
\item $t\in I$, $i,k\not\in I$, $m_{ik}=1$, then with
$*)$: cases $(a)$ and $(c)$ of $ii)$ above
$$
\tilde b_{ik} \tilde b_t \stackrel{*)}{\in}
\tilde b_{ik} \tilde b_i \tilde b_k\inv \tilde b_t H_M \stackrel{iv)}{=}
\tilde b_k \tilde b_{ik} \tilde b_k\inv \tilde b_t H_M  \stackrel{ii)}{=}
\tilde b_t \tilde b_k \tilde b_{ik} \tilde b_k\inv H_M \stackrel{}{=}
\tilde b_t H_M
$$
\item $t\,i,k\in I$, i.e. $m_{ik}=1$, then with
$*)$: cases $(a)$ and $(c)$ of $i)$ above
$$
\tilde b_{ik} \tilde b_t \stackrel{*)}{\in}
\tilde b_{ik} \tilde b_i^2 \tilde b_k\invv \tilde b_t H_M \stackrel{v)}{=}
\tilde b_k^2 \tilde b_{ik} \tilde b_k\invv \tilde b_t H_M \stackrel{ii)}{=}
\tilde b_t \tilde b_k^2 \tilde b_{ik} \tilde b_k\invv H_M \stackrel{}{=}
\tilde b_t H_M
$$
\end{enumerate}
\end{enumerate}
\qed

We are now ready to give the proof of our second theorem.
\\

\proof[ of thm.\ \ref{manf}]
The case $n=1$ is void and the case $n=2$ consists of two subcases
with group $E$ freely generated by $\s_1$ resp.\ $\s_1^2$ in
accordance with the claim.
Hence it suffices to prove the claim for matrices $M$ of size $n$
relying on the induction hypothesis. In fact we will show that there
is an isomorphism $G_M\to E_M$ induced by the natural bijection
$\tilde b_{ij}\to b_{ij}$ on generators.

First we observe that we get in fact a surjection, since the relations
in $G_M$ map to relations in $\br_n$, which can be checked
by straightforward calculations.
Moreover we note that the subgroup $H_M$ maps onto
$E_{P'n}$ where $P'n$ is obtained as before from the partition $P$
associated to $M$. 

We will now exploit the induction hypothesis and
lemma \ref{combing}, which provide a presentation of $E_{P'n}$, to
get at least a partial inverse $E_{P'n}\to H_M$.
Since the relations in $\rfami'$ pose no problems we restrict to
enumerate the relations in $G_M$ from which the 
relations of $E_{P'n}$ in the corresponding
enumeration of lemma \ref{combing} follow.
\begin{enumerate}
\item
$\tilde b_i \tilde b_{jk} = \tilde b_{jk} \tilde b_{i}$
\item according to $m_{kn}=2$ or $1$ either $a_k^2\mapsto \tilde b_{k}$ (a),
or $a_k^2\mapsto \tilde b_{k}^2$ (b).
\begin{enumerate}
\item 
$\tilde b_j \tilde b_k \tilde b_{ik} \tilde b_k\inv
= \tilde b_k \tilde b_{ik} \tilde b_k\inv \tilde b_j$
\item 
$\tilde b_j \tilde b_k^2 \tilde b_{ik} \tilde b_k\invv
= \tilde b_k^2 \tilde b_{ik} \tilde b_k\invv \tilde b_j$
\end{enumerate}
\item since $m_{ij}=1$ either $a_i,a_j\not\in E_M$
and $a_i^2\mapsto \tilde b_i,\,a_j^2\mapsto \tilde b_j$ (a),\\ 
or $a_i,a_j\in E_M$ and $a_i^2\mapsto \tilde b_i^2$, 
$a_j^2\mapsto \tilde b_j^2$ (b).
\begin{enumerate}
\item
$\tilde b_{ij} \tilde b_{j} = \tilde b_j \tilde b_{ij},\,
\tilde b_{ij} \tilde b_i \tilde b_j =  \tilde b_i \tilde b_j \tilde b_{ij}$
\item
$\tilde b_{ij} \tilde b_i = \tilde b_i \tilde b_j
= \tilde b_j \tilde b_{ij}$ implies
$\tilde b_{ij} \tilde b_i^2 = \tilde b_j^2 \tilde b_{ij}$ and\\
$\tilde b_{ij} \tilde b_i^2 \tilde b_j^2 
= \ul{\tilde b_{ij} \tilde b_i}\,\ul{ \tilde b_i \tilde b_j} \tilde b_j 
= \tilde b_i \ul{\tilde b_j \tilde b_{ij}}\,\ul{ \tilde b_i \tilde b_j} 
= \tilde b_i \tilde b_i \tilde b_j \tilde b_j \tilde b_{ij} 
= \tilde b_i^2 \tilde b_j^2 \tilde b_{ij}$
\end{enumerate}
\item since $m_{ij}=2$ either $a_i,a_j\not\in E_M$
and $a_i^2\mapsto \tilde b_i,\,a_j^2\mapsto \tilde b_j$ (a),\\ 
or $a_i\in E_M, a_j\not \in E_M$ 
and $a_i^2\mapsto \tilde b_i^2$, $a_j^2\mapsto \tilde b_j^2$ (b), \\
or $a_i\not\in E_M, a_j \in E_M$ 
and $a_i^2\mapsto \tilde b_i^2$, $a_j^2\mapsto \tilde b_j^2$ (c).
\begin{enumerate}
\item 
$\tilde b_i \tilde b_j \tilde b_{ij} 
= \tilde b_j \tilde b_{ij} \tilde b_i$,\,
$\tilde b_i \tilde b_j \tilde b_{ij} 
= \tilde b_{ij} \tilde b_i \tilde b_j$
\item 
$\tilde b_i \tilde b_j \tilde b_{ij} 
= \tilde b_j \tilde b_{ij} \tilde b_i$
and $\tilde b_i \tilde b_j = \tilde b_{ij} \tilde b_i$ imply\\
$\tilde b_i^2 \tilde b_j \tilde b_{ij} 
= \tilde b_j \tilde b_{ij} \tilde b_i^2$ and
$\tilde b_i^2 \tilde b_j \tilde b_{ij} 
= \ul{\tilde b_i \tilde b_j}\, \ul{\tilde b_{ij} \tilde b_i} 
= \tilde b_{ij} \tilde b_i \tilde b_i \tilde b_j 
= \tilde b_{ij} \tilde b_i^2 \tilde b_j$
\item
$\tilde b_j \tilde b_{ij} \tilde b_i 
= \tilde b_{ij} \tilde b_i \tilde b_j$ and 
$\tilde b_j \tilde b_{ij} = \tilde b_i \tilde b_j$ imply\\
$\tilde b_j^2 \tilde b_{ij} \tilde b_i 
= \tilde b_{ij} \tilde b_i \tilde b_j^2$ and
$\tilde b_j^2 \tilde b_{ij} \tilde b_i 
= \ul{\tilde b_j \tilde b_{ij}}\, \ul{\tilde b_i \tilde b_j} 
= \tilde b_i \tilde b_j \tilde b_j \tilde b_{ij} 
= \tilde b_i \tilde b_j^2 \tilde b_{ij}$
\end{enumerate}
\end{enumerate}

Next we consider the following diagram with exact rows:
$$
\begin{array}{ccccccccc}
1 & \to & \pbr_n & \to & E_{P'n} & \to & S_{n,P'n} & \to & 1\\
&& \| && \cap && \cap\\
1 & \to & \pbr_n & \to & E_P & \to & S_{n,P} & \to & 1
\end{array}
$$
Hence $E_{P'n}$ is of finite index 
$|I|=\#\{i\leq n\,|\,i\sim_P n\}$ in $E_P$.
This injection factors as 
$E_{P'n}\cong H_M\inj G_M \surj E_P$.
Since by lemma \ref{cosets} the index of $H_M$ in $G_M$ is at
most $|I|$, we can conclude that $E_P$ is isomorphic to $G_M$.
\qed

\section{further prospects}
\labell{unfolding}

In this last section we want to address two directions
of further studies.
First we would like to mention the simplest instance of a structural
result describing groups associated to a matrix $M$
composed of several submatrices in terms of the
groups associated to these submatrices.

\begin{prop}
Suppose $M$ is a block matrix with diagonal blocks $M_1, M_2$
and zero off-diagonal blocks.
Then the groups associated to $M_1,M_2$ are free or direct factors:
$$
A_M = A_{M_1}*A_{M_2},\quad
C_M = C_{M_1}*C_{M_2},\quad
E_M = E_{M_1}\times E_{M_2}
$$
\end{prop}

It should also be worthwhile to understand how the group $E_M$
changes under the action of a permutation matrix on $M$.
\\

More emphasis we would like to put on the problem to find
presentations for matrices with entries $m_{ij}\geq 2$, in
particular with matrices of finite irreducible Coxeter groups.

The following result on relation between the generators of
such groups $E_M$ is obtained by straightforward calculations.

\begin{prop}
In case of matrices $M$ with entries $m_{ij}\geq2$, 
the group $E_M$ generated by $b_{ij}=a_{ij}^{m_{ij}}$
has relations
\begin{enumerate}
\item
for $i<j<k<l$
$$
b_{ij} b_{kl} = b_{kl} b_{ij},\quad
b_{il} b_{jk} = b_{jk} b_{il},
$$
\item
for $i<j<k<l$ and $m_{jk}=2$
$$
b_{ik} b_{jk} b_{jl} b_{jk}\inv = b_{jk} b_{jl} b_{jk}\inv b_{ik}.
$$
\item
for $i<j<k$, $j<k<i$ or $k<i<j$ with
\begin{enumerate}
\item 
$m_{ij}=m_{ik}=2$, $m_{jk}=2\nu$
$$
(b_{ij} b_{ik})^{\nu-1} b_{jk} b_{ij} b_{ik} = 
b_{ik} (b_{ij} b_{ik})^{\nu-1} b_{jk} b_{ij} = 
(b_{ij} b_{ik})^\nu b_{jk},
$$
\item 
$m_{ij}=m_{ik}=2$, $m_{jk}=2\nu+1$
$$
b_{ik} (b_{ij} b_{ik})^{\nu-1} b_{jk} b_{ij} b_{ik} 
= (b_{ij} b_{ik})^{\nu-1} b_{jk} b_{ij} 
= b_{ik} (b_{ij} b_{ik})^{\nu} b_{jk},
$$

\item
$m_{ij}=2, m_{ik}=m_{jk}=3$
$$
b_{ij} b_{jk} b_{ij} b_{ik} b_{jk} = 
b_{jk} b_{ij} b_{ik} b_{jk} b_{ij} = b_{ij} b_{ik} b_{jk} b_{ij} b_{ik},
$$

\item
$m_{ij}=2,m_{ik}=3,m_{jk}=4$
$$
b_{ij} b_{ik} b_{jk} b_{ij} b_{ik} b_{jk} =
b_{ik} b_{jk} b_{ij} b_{ik} b_{jk} b_{ij} = b_{jk} b_{ij} b_{ik} b_{jk} b_{ij} b_{ik},
$$
\item
$m_{ij}=2,m_{ik}=3,m_{jk}=5$
$$
b_{ij} b_{ik} b_{jk} b_{ij} b_{ik} b_{jk} b_{ij} b_{ik} b_{jk} =
b_{ik} b_{jk} b_{ij} b_{ik} b_{jk} b_{ij} b_{ik} b_{jk} b_{ij} = 
b_{jk} b_{ij} b_{ik} b_{jk} b_{ij} b_{ik} b_{jk} b_{ij} b_{ik}.
$$
\end{enumerate}

\end{enumerate}
\end{prop}

In fact we will prove that the given relations are sufficient to present the
group $E_M$ in case $M$ is a Coxeter matrix for a Coxeter group of type $A$,
$D$ or $E$ in a forthcoming paper.
The methods employed will be totally different, since we will actually exploit
braid monodromy techniques, cf.\ \cite{habil}, to plane curve complements, which
can be shown to have complements with fundamental group isomorphic
to $E_M$.


\begin{thebibliography}{XXX}


\bibitem[Ar]{ar} Emil Artin: {\sl Theory of Braids},
Ann. Math. (2) 48 (1947), 101--126

\bibitem[Bu]{bu} Walter Burau: {\sl \"Uber Zopfinvarianten},
Abh. Math. Semin. Hamb. Univ. 9 (1932), 117--124

\bibitem[CP]{cp} John Crisp, Luis Paris:
{\sl The solution to a conjecture of Tits on the subgroup generated by the squares of the generators of an Artin group}, 
Invent. Math.  145  (2001),  no. 1, 19--36. 

\bibitem[D\"o]{Doe}
Axel D\"orner:
{\sl Isotropieuntergruppen der artinschen Zopfgruppen},
Dissertation; Bonner Mathematische Schriften 255,
Mathematisches Institut, Bonn (1993)

\bibitem[Kl]{klui} Paul Kluitmann:
{\sl Isotropy Subgroups of Artin's Braid
Group}, preprint 1991

\bibitem[L\"o1]{hurwitz} Michael L\"onne: {\sl Hurwitz stabilizers of some short
redundant Artin systems for the braid group $\br_3$}, preprint AG/0406154

\bibitem[L\"o2]{habil} Michael L\"onne: {\sl Braid monodromy of 
hypersurface singularities},
Habilitations\-schrift (2003), Hannover, mathAG/0602371,

\bibitem[Lo]{loo} Eduard Looijenga: {\sl The complement of the bifurcation variety
of a simple singularity},
Invent.\ Math.\ 23 (1974), 105--116

\bibitem[Man]{manf} Sandro Manfredini: {\sl Some subgroups of Artin's braid group},
Topology Appl. 78 (1997), 123--142


\end{thebibliography}
\end{document}